%% file: MIFBM_160615.tex
%
\newif\ifloadreferences\loadreferencestrue
%
\input preamble %
%
%
%
%
%
\def\Pagetitle{\hfil}
\def\Pagefooter{\hfil{\myfontdefault\folio}\hfil}
\null \vfill
\def\centre{\rightskip=0pt plus 1fil \leftskip=0pt plus 1fil \spaceskip=.3333em \xspaceskip=.5em \parfillskip=0em \parindent=0em}%
\def\textmonth#1{\ifcase#1\or January\or Febuary\or March\or April\or May\or June\or July\or August\or September\or October\or November\or December\fi}
\font\abstracttitlefont=cmr10 at 14pt {\abstracttitlefont\centre
The Morse index of the critical catenoid.\par}
\bigskip
{\centre 7th September 2016\par}
\bigskip
{\centre Graham Smith\par}
\medskip
{\centre Instituto de Matem\'atica,\par
UFRJ, Av. Athos da Silveira Ramos 149,\par
Centro de Tecnologia - Bloco C,\par
Cidade Universit\'aria - Ilha do Fund\~ao,\par
Caixa Postal 68530, 21941-909,\par
Rio de Janeiro, RJ - BRASIL\par}
\medskip
{\centre Detang Zhou\par}
\medskip
{\centre Instituto de Matem\'atica,\par
UFF, Rua Professor Marcos Waldemar de Freitas Reis,\par
Bloco H - Campus do Gragoat\'a,\par
S\~ao Domingos, 24.210-201,\par
Niter\'oi, RJ - BRASIL\par}
\bigskip
\noindent{\bf Abstract:\ }We show that the rotationally symmetric free boundary minimal catenoid in the unit ball in $\Bbb{R}^3$ has Morse index equal to $4$.
\par
\vfill
\eject
%
\myfontdefault
\def\Pagetitle{\hfil The Morse index of the critical catenoid.\hfil}
\def\Pagefooter{\hfil{\myfontdefault\folio}\hfil}
\makeop{rob}%
\makeopsmall{Coth}%
\makeop{dA}%
\makeop{dl}%
\newref{Brendle}{Brendle S., Embedded minimal tori in $S^3$ and the Lawson conjecture, {\sl Acta Math.}, {\bf 211}, no. 2, (2013), 177--190}
\newref{ChenFraserPang}{Chen J., Fraser A., Pang C., Minimal immersions of compact bordered Riemann surfaces with free boundary, {\sl Trans. Amer. Math. Soc.}, {\bf 367}, (2015), 2487--2507}
\newref{CodaNeves}{Cod\'a Marques F., Neves A., Min-max theory and the Willmore conjecture, {\sl Ann. Math.}, {\bf 179}, no. 2, (2014), 683--782}
\newref{Devyver}{Devyver B., Index of the critical catenoid, arXiv:1609.02315}
\newref{FraserSchoenI}{Fraser A., Schoen R., The first Steklov eigenvalue, conformal geometry, and minimal surfaces., {\sl Adv. Math.}, {\bf 226}, (2011), no. 5, 4011--4030}
\newref{FraserSchoenII}{Fraser A., Schoen R., Sharp eigenvalue bounds and minimal surfaces in the ball, {\sl Invent. Math.}, {\bf 203}, (2016), no. 3, 823--890}
\newref{MaximoNunesSmith}{M\'aximo D., Nunes I., Smith G., Free boundary minimal annuli in convex three-manifolds, to appear in {\sl J. Diff. Geom.}}
\newref{Hung}{Tranh H. T., Index characterization for free boundary minimal surfaces, arXiv:1609.01651}
\newref{Urbano}{Urbano F., Minimal surfaces with low index in the three-dimensional sphere, {\sl Proc.
Amer. Math. Soc.}, {\bf 108}, (1990), 989--992}
%
%
\newsubhead{The Morse index of the critical catenoid}[MorseIndexOfTheCriticalCatenoid]
Let $B:=B^3$ be the closed unit ball in $3$-dimensional Euclidean space. Let $\Sigma$ be a compact surface smoothly embedded in $B$ with smooth boundary in $\partial B$. We say that $\Sigma$ is {\sl free boundary minimal} whenever it is a critical point of the area functional amongst all such embedded surfaces in $B$. It is well known that this implies that $\Sigma$ has vanishing mean curvature, that $\Sigma$ meets $\partial B$ orthogonally along $\partial\Sigma$, and that no interior point of $\Sigma$ meets $\partial B$. Free boundary minimal surfaces have been studied since the classical work of Gergonne, and have received renewed attention with the recent work \cite{FraserSchoenI} and \cite{FraserSchoenII} of Fraser \& Schoen.
\par
Infinitesimal perturbations of $\Sigma$ are given by normal vector fields. If we suppose that $\Sigma$ is orientable, then we choose a unit normal vector field $N:\Sigma\rightarrow S^2$ which is compatible with the orientation. Using this vector field, we now identify infinitesimal perturbations of $\Sigma$ with smooth functions over this surface by identifying the field $fN$ with the function $f$. Then, given $f$, the second order variation of area for an infinitesimal perturbation in the direction of $f$ is given by (c.f. \cite{MaximoNunesSmith})
$$
S(f,f) := \int_\Sigma \|\nabla f\|^2 - \opTr(A)^2f^2\opdA - \int_{\partial\Sigma}f^2\opdl,
$$
where $\nabla$ denotes the gradient operator of $\Sigma$, $A$ its shape operator, $\opdA$ its area element and $\opdl$ the length element of $\partial\Sigma$. It is of interest to study the spectrum of this symmetric bilinear form, where a function $f:\Sigma\rightarrow\Bbb{R}$ is defined to be an eigenfunction of $S$ with eigenvalue $\lambda$ whenever
$$
S(f,g) = \lambda\langle f,g\rangle
$$
for any other function $g$, where $\langle\cdot,\cdot\rangle$ here denotes the $L^2$ inner product of $\Sigma$. However, it is straightforward to show (c.f. \cite{ChenFraserPang}, Section $5$) that this spectrum coincides with the spectrum of the operator
$$
Jf := -\Delta f - \opTr(A^2)f,
$$
restricted to the space of functions $f$ which satisfy the {\sl Robin boundary condition}
$$
f=\partial_\nu f,\eqnum{\nexteqnno[RobinBoundaryCondition]}
$$
where here $\Delta$ denotes the laplacian operator of $\Sigma$ and $\nu$ denotes the outward pointing unit conormal vector field over $\partial\Sigma$. The infinitesimal perturbations of $\Sigma$ given by functions which satisfy \eqnref{RobinBoundaryCondition} have a straightforward geometric interpretation: they are precisely those perturbations which preserve up to first order the property of $\Sigma$ meeting $\partial B$ orthogonally.
\par
In this note, we will be interested in the case where $\Sigma$ is the critical catenoid, which we denote by $\Sigma_c$. Recall that this is the surface of revolution about the $x$-axis of the function
$$
\frac{1}{R}\opCosh(Rx),
$$
over the interval
$$
[-L/R,L/R],
$$
where
$$
R:=L\cdot\opCosh(L),
$$
and $L$ is the unique positive solution of
$$
L=\opCoth(L).
$$
Of particular interest are the {\sl Morse index} and {\sl nullity} of $\Sigma_c$, which are defined respectively to be the number of strictly negative eigenvalues of $S$ - counted with multiplicity - and the dimension of its null space. In \cite{MaximoNunesSmith}, we showed that the nullity of $\Sigma_c$ is equal to $2$. In the present note, we show
\proclaim{Theorem \nextprocno}
\noindent The Morse index of $\Sigma_c$ is equal to $4$.
\endproclaim
\proclabel{MorseIndex}
Although the study of embedded, free boundary minimal catenoids in $B^3$ would seem to be analogous to the study of embedded minimal tori in the $3$-sphere $S^3$, it is actually much harder. Indeed, wheras long-standing conjectures concerning tori in $S^3$ have recently been resolved, the corresponding problems in the free boundary case remain completely open. For example, in \cite{Brendle}, Brendle solves the Lawson conjecture by showing that the so called Clifford tori are, in fact, the only embedded minimal tori in $S^3$. However, his argument breaks down completely in the free boundary case, and no substitute has yet been found. Likewise, in \cite{CodaNeves}, Coda \& Neves solve the Wilmore conjecture, proving in the process that the Clifford tori minimise area amongst all embedded minimal submanifolds of $S^3$ that are not topologically spheres. However, their argument makes use of the fact that the Wilmore energy of an embedded surface in $S^3$ is bounded below by its area, making it unclear how it should adapt to the free boundary case, where no such relationship holds.
\par
Of particular relevance to the present note is the result \cite{Urbano} of Urbano, which Coda \& Neves showed to be of considerable use in their proof of the Wilmore conjecture. Recall that if $S^3$ is identified with the unit sphere in $4$-dimensional Euclidean space, that is
$$
S^3 = \left\{ x_1^2 + x_2^2 + x_3^2 + x_4^2 = 1\right\},
$$
then the Clifford tori are the images under the action of the rotation group of the surface
$$
C := \left\{ x_1^2 + x_2^2 = x_3^2 + x_4^2 = \frac{1}{2}\right\}.
$$
It is a straightforward exercise to show that the Morse index of a Clifford torus is equal to $5$. Urbano proves, conversely, that this property characterises Clifford tori amongst embedded minimal surfaces in $S^3$. Indeed, if $\Sigma$ is an embedded minimal surface in $S^3$ of Morse index less than or equal to $5$, then it is either an equatorial sphere or a Clifford torus.
\par
One would conjecture a similar result for embedded, free boundary minimal annuli in $B^3$. However, it seems unlikely that Urbano's argument would apply as, like the work of Coda \& Neves, his proof makes non-trivial use of the Wilmore energy which, as we have already indicated, is unrelated to the area in the free boundary case. In fact, we are not aware that anyone has even determined the Morse index of the critical catenoid, hence the current note.
\medskip
Upon completion of this paper, we were made aware of similar results \cite{Devyver} by Baptiste Devyver and \cite{Hung} by Hung Tran obtained simultaneously with our own using completely different techniques.
\newsubhead{Proof of Theorem \procref{MorseIndex}}[ProofOfTheoremX]
This is a more or less straightforward application of the technique of separation of variables. We study the spectrum of $J$ over the space of functions satisfying the Robin boundary condition \eqnref{RobinBoundaryCondition}. In particular, by elliptic regularity, all eigenfunctions are smooth up to the boundary. Consider now the parametrisation $\Phi:S^1\times[-L,L]\rightarrow\Sigma_c$ given by
$$
\Phi(\theta,x) = \frac{1}{R}(x,\opCosh(x)\opCos(\theta),\opCosh(x)\opSin(\theta)).
$$
It is a straightforward exercise (c.f. \cite{MaximoNunesSmith}) to show that in these coordinates, a function $f:S^1\times[-L,L]\rightarrow\Bbb{R}$ is an eigenfunction of $J$ with eigenvalue $-\lambda$ whenever
$$
\frac{1}{\opCosh^2(x)}(f_{xx} + f_{\theta\theta}) + \frac{2}{\opCosh^4(x)}f = \lambda f,\eqnum{\nexteqnno[EigenfunctionEquation]}
$$
and $f$ satisfies the Robin boundary condition \eqnref{RobinBoundaryCondition} whenever
$$\eqalign{
f(L) &= Lf'(L),\ \text{and}\cr
f(-L) &= -Lf'(-L).\cr}\eqnum{\nexteqnno[OneDimRobin]}
$$
We aim to determine for which positive values of $\lambda$ this problem has a non-trivial solution. To this end, we introduce the ansatz,
$$
f(\theta,x) = \sum_{m\in\Bbb{Z}}e^{im\theta}f_m(x).
$$
The function $f$ then solves \eqnref{EigenfunctionEquation} and \eqnref{OneDimRobin} whenever, for every integer, $m$,
$$
f''_m = \left(m^2 + \lambda\opCosh^2(x) - \frac{2}{\opCosh^2(x)}\right)f_m,\eqnum{\nexteqnno[EEInFourierMode]}
$$
and
$$\eqalign{
f_m(L) &= Lf_m'(L),\ \text{and}\cr
f_m(-L) &= -Lf_m'(-L).\cr}\eqnum{\nexteqnno[RCInFourierMode]}
$$
Significantly, if $f_m(x)$ solves this problem, then so too does $f_m(-x)$, so that, given any solution $f_m$ its even and odd components are also solutions.
\proclaim{Lemma \nextprocno}
\noindent If $m^2+\lambda>3$, then there are no non-trivial solutions to \eqnref{EEInFourierMode} and \eqnref{RCInFourierMode}.
\endproclaim
\proclabel{UpperLimits}
\proof Let $f$ solve the differential equation \eqnref{EEInFourierMode}, and suppose that $f$ is either even or odd. We will show that $f$ cannot satisfy the Robin boundary conditions \eqnref{RCInFourierMode}. Indeed, consider the function
$$
\gamma := \frac{f'}{f}.
$$
This function satisfies the Ricatti equation
$$
\gamma' + \gamma^2 = \left( m^2 + \lambda\opCosh^2(x) - \frac{2}{\opCosh^2(x)}\right) > 1.\eqnum{\nexteqnno[RicattiEquation]}
$$
There are now two cases to consider. If $f$ is even, then $\gamma(0)=0$. By \eqnref{RicattiEquation}, $\gamma$ can only become singular by becoming large and negative. However, this relation also implies that
$$
\gamma(x) > \opTanh(x).
$$
In particular, $\gamma$ is never singular, and $\gamma(L)>\opTanh(L)=L^{-1}$, so that $f$ does not satisfy the Robin boundary condition in this case.
\par
If $f$ is odd, then, without loss of generality,
$$
\mlim_{x\rightarrow 0}\left(\gamma(x) - \frac{1}{x}\right)= 0.
$$
As before, it follows by \eqnref{RicattiEquation} that $\gamma$ is non-singular for all $x>0$, and that
$$
\gamma(x) > \opCoth(x).
$$
In particular, $\gamma(L)>\opCoth(L)=L>L^{-1}$, and so $f$ does not satisfy the Robin boundary condition in this case either. This completes the proof.\qed
\medskip
It follows by Lemma \procref{UpperLimits} that eigenfunctions of $J$ only have finitely many non-trivial Fourier modes, and eigenfunctions with negative eigenvalues only have non-trivial Fourier modes of orders $0$ or $\pm 1$.
\par
The low order Fourier modes are now studied in terms of the geometry of $\Sigma$. Consider first the $0$'th order mode, and define
$$\eqalign{
f(x) &:= 1 - x\opTanh(x),\ \text{and}\cr
g(x) &:= \opTanh(x).\cr}
$$
These functions are even and odd respectively, and both solve \eqnref{EEInFourierMode} with $m=0$ and $\lambda=0$. Geometrically, $f$ and $g$ are respectively the infinitesimal perturbations of $\Sigma_c$ given by dilatations centred on the origin, and translations along the $x$-axis. We now define a family of functions parametrised by $[a:b]\in\Bbb{R}\Bbb{P}^1$ by
$$
\gamma_{[a:b]}(x) := \frac{af'(x) + bg'(x)}{af(x) + bg(x)}.
$$
\placefigure{}{}{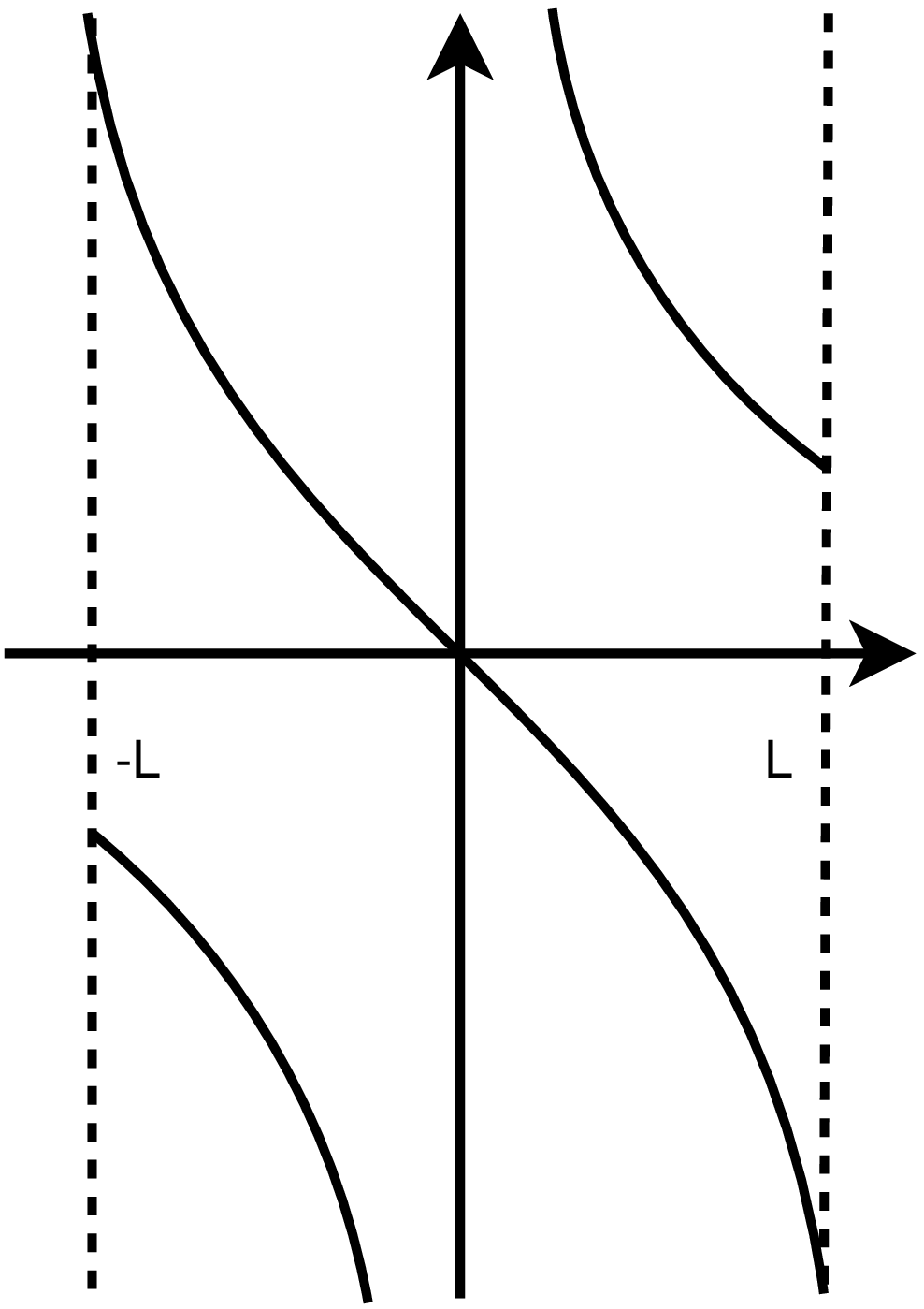}{The rectangle $[-L,L]\times\Bbb{R}$ is foliated by integral curves of the vector field $X$. Shown here are the two leaves defined by the graphs of $\gamma_{[1:0]}$, passing through the origin, and $\gamma_{[0:1]}$, asymptotic to the $y$-axis.}
\figlabel{FigureA}
The graphs of these functions define a smooth foliation $\Cal{F}$ of the rectangle $[-L,L]\times\Bbb{R}$ whose leaves are parametrised by $\Bbb{R}\Bbb{P}^1$ (c.f. Figure \figref{FigureA}). In particular, this foliation integrates the vector field
$$
X(x,y) := \left(1, -\frac{2}{\opCosh^2(x)} - y^2\right).
$$
\par
Consider now the function $\gamma_{[1:0]}$. This corresponds to the even solution $f$ of \eqnref{EEInFourierMode} and satisfies
$$\eqalign{
\gamma_{[1:0]}(x) &= 0,\ \text{and}\cr
\gamma_{[1:0]}(x) &\rightarrow-\infty\ \text{as}\ x\rightarrow L.\cr}
$$
For $\lambda>0$, let $\gamma_{[1:0],\lambda}$ be the solution of the Ricatti equation
$$
\gamma'(x) + \gamma(x)^2 =\lambda\opCosh^2(x) - \frac{2}{\opCosh^2(x)},
$$
with initial condition $\gamma(0)=0$. Since the graph of this function is an integral curve of the vector field
$$
X_\lambda(x,y) := \left(1,\lambda\opCosh^2(x) - \frac{2}{\opCosh^2(x)} - y^2\right),
$$
it only ever crosses the foliation $\Cal{F}$ in the upward direction. Thus, since it can only become singular by becoming large and negative, it follows that $\gamma_{[1:0],\lambda}$ is finite and smooth over the interval $[0,L]$ for all $\lambda>0$. In particular, if we define the function $\phi:]0,\infty[\rightarrow\Bbb{R}$, by
$$
\phi(\lambda) = \gamma_{[1:0],\lambda}(L),
$$
then $\phi$ is a strictly increasing, smooth function such that
$$
\phi(\lambda)\rightarrow-\infty\ \text{as}\ \lambda\rightarrow 0.
$$
However, by the proof of Lemma \procref{UpperLimits},
$$
\phi(3)>\frac{1}{L},
$$
and it follows by the intermediate value theorem that there exists a unique, positive value of $\lambda$ such that
$$
\phi(\lambda) = \gamma_{[1:0],\lambda}(L) = \frac{1}{L}.
$$
This $\lambda$ is an eigenvalue of $J$ whose eigenfunction is even and lies in the $0$'th Fourier mode.
\par
Consider now the function, $\gamma_{[0:1]}$. This corresponds to the odd solution $g$ of \eqnref{EEInFourierMode} and satisfies
$$\eqalign{
\mlim_{x\rightarrow 0}\left(\gamma_{[0:1]}(x) - \frac{1}{x}\right) &=0,\ \text{and}\cr
\gamma_{[0:1]}(L) &= L - \frac{1}{L}.\cr}
$$
We define the functions $\gamma_{[0:1],\lambda}$ and $\phi$ as before, and we see that $\phi$ is now a strictly increasing, smooth function over $[0,\infty[$ such that
$$
\phi(0) = L - \frac{1}{L} < \frac{1}{L}.
$$
Since $\phi(3)>L^{-1}$, it follows again by the intermediate value theorem that there exists a unique, positive value of $\lambda$ such that
$$
\phi(\lambda) = \gamma_{[0:1],\lambda}(L) = \frac{1}{L}.
$$
This $\lambda$ is an eigenvalue of $J$ whose eigenfunction is odd and lies in the $0$'th Fourier mode.
\par
The Fourier modes of order $\pm 1$ are treated in a similar manner. Indeed, define
$$\eqalign{
f(x) &:= \frac{1}{\opCosh(x)},\ \text{and}\cr
g(x) &:= \opSinh(x) + \frac{x}{\opCosh(x)}.\cr}
$$
These functions are even and odd respectively, and both solve \eqnref{EEInFourierMode} with $m=1$ and $\lambda=0$. Geometrically, $f$ and $g$ are respectively the infinitesimal perturbations given by translations in directions orthogonal to the $x$-axis, and rotations about axes orthogonal to the $x$-axis.
\par
Proceeding as before, for all $[a:b]\in\Bbb{R}\Bbb{P}^1$, we define the function $\gamma_{[a:b]}$. In this case, $\gamma_{[1:0]}(L)=-L^{-1}<L^{-1}$, and thus yields a unique negative eigenvalue with two eigenfunctions that are even and lie in the Fourier modes of order $-1$ and $1$. However, $\gamma_{[0:1]}(L)=L^{-1}$, so that the perturbations of this solution yield no new negative eigenvalues.
\par
In summary, we have shown that $J$ has $2$ negative eigenvalues in the $0$'th Fourier mode, and $1$ negative eigenvalue in each of the Fourier modes of order $-1$ and $1$. This proves Theorem \procref{MorseIndex}.
\newsubhead{Bibliography}[Bibliography]
\medskip
{\leftskip = 5ex \parindent = -5ex
\leavevmode\hbox to 4ex{\hfil \cite{Brendle}}\hskip 1ex{Brendle S., Embedded minimal tori in $S^3$ and the Lawson conjecture, {\sl Acta Math.}, {\bf 211}, no. 2, (2013), 177--190}
\medskip
\leavevmode\hbox to 4ex{\hfil \cite{ChenFraserPang}}\hskip 1ex{Chen J., Fraser A., Pang C., Minimal immersions of compact bordered Riemann surfaces with free boundary, {\sl Trans. Amer. Math. Soc.}, {\bf 367}, (2015), 2487--2507}
\medskip
\leavevmode\hbox to 4ex{\hfil \cite{CodaNeves}}\hskip 1ex{Marques F. C., Neves A., Min-max theory and the Willmore conjecture, {\sl Ann. Math.}, {\bf 179}, no. 2, (2014), 683--782}
\medskip
\leavevmode\hbox to 4ex{\hfil \cite{Devyver}}\hskip 1ex{Devyver B., Index of the critical catenoid, arXiv:1609.02315}
\medskip
\leavevmode\hbox to 4ex{\hfil \cite{FraserSchoenI}}\hskip 1ex{Fraser A., Schoen R., The first Steklov eigenvalue, conformal geometry, and minimal surfaces., {\sl Adv. Math.}, {\bf 226}, (2011), no. 5, 4011--4030}
\medskip
\leavevmode\hbox to 4ex{\hfil \cite{FraserSchoenII}}\hskip 1ex{Fraser A., Schoen R., Sharp eigenvalue bounds and minimal surfaces in the ball, {\sl Invent. Math.}, {\bf 203}, (2016), no. 3, 823--890}
\medskip
\leavevmode\hbox to 4ex{\hfil \cite{MaximoNunesSmith}}\hskip 1ex{M\'aximo D., Nunes I., Smith G., Free boundary minimal annuli in convex three-manifolds, to appear in {\sl J. Diff. Geom.}}
\medskip
\leavevmode\hbox to 4ex{\hfil \cite{Hung}}\hskip 1ex{Tran H., Index characterization for free boundary minimal surfaces, arXiv:1609.01651}
\medskip
\leavevmode\hbox to 4ex{\hfil \cite{Urbano}}\hskip 1ex{Urbano F., Minimal surfaces with low index in the three-dimensional sphere, {\sl Proc. Amer. Math. Soc.}, {\bf 108}, (1990), 989--992}
\par}
%
%
%
%
\enddocument

%% file: preamble.tex
%
%
%
%
\let\myfrac=\frac%
\input eplain %
\let\frac=\myfrac%
\let\myfootnote=\footnote%
\input amstex \input epsf %
\let\footnote=\myfootnote%
%
%
\loadeufm\loadmsam\loadmsbm\message{symbol names}\UseAMSsymbols\message{,}%
\magnification 1200 %
\font\myfontdefault=cmr10%
\newif\ifmakebiblio%
\newif\ifinappendices%
\newif\ifundefinedreferences%
\newif\ifchangedreferences%
\makebibliofalse%
\undefinedreferencesfalse%
\changedreferencesfalse%
%
%
%
%
%
\def\setcatcodes{\catcode`\!=0 \catcode`\\=11}%
{\global\let\noe=\noexpand%
\catcode`\@=11 \catcode`\_=11 \setcatcodes%
!global!def!_@@internal@@makeref#1{%
!global!expandafter!def!csname #1ref!endcsname##1{%
!csname _@#1@##1!endcsname%
!expandafter!ifx!csname _@#1@##1!endcsname!relax%
    !write16{#1 ##1 not defined - run saving references}%
    !undefinedreferencestrue%
!fi}}%
!global!def!_@@internal@@makelabel#1{%
!global!expandafter!def!csname #1label!endcsname##1{%
!edef!temptoken{!csname #1info!endcsname}%
!ifloadreferences%
!expandafter!ifx!csname _@#1@##1!endcsname!relax%
!write16{#1 ##1 not hitherto defined - rerun saving references}%
!changedreferencestrue%
!else%
!expandafter!ifx!csname _@#1@##1!endcsname!temptoken%
!else%
!write16{#1 ##1 reference has changed - rerun saving references}%
!changedreferencestrue%
!fi%
!fi%
!else%
!expandafter!edef!csname _@#1@##1!endcsname{!temptoken}%
!edef!textoutput{!write!references{\global\def\_@#1@##1{!temptoken}}}%
!textoutput%
!fi}}%
!global!def!makecounter#1{!_@@internal@@makelabel{#1}!_@@internal@@makeref{#1}}%
!unsetcatcodes%
}
%
%
%
%
%
\def\turnintolatin#1{\ifcase #1 _\or i\or ii\or iii\or iv\or v\or vi\or vii\or viii\or ix\or x\or xi\or xii\or xiii\or xiv\or xv\or xvi\or xvii\or xviii\or xix\or xx\or xxi\or xxii\or xxiii\or xxiv\or xxv\or xxvi\fi}%
\def\alphanum#1{\ifcase #1 _\or A\or B\or C\or D\or E\or F\or G\or H\or I\or J\or K\or L\or M\or N\or O\or P\or Q\or R\or S\or T\or U\or V\or W\or X\or Y\or Z\fi}%
\newwrite\references%
\ifloadreferences{\catcode`\@=11 \catcode`\_=11 \input references.tex }%
\else{\openout\references=references.tex }%
\fi%
%
%
\newcount\headno%
\global\headno=0%
\def\headinfo{\ifinappendices\alphanum\headno\else\the\headno\fi}%
\def\nextheadno{\global\advance\headno by 1 \global\subheadno=0 \global\procno=0 \global\eqnno=0 \headinfo}%
\makecounter{head}%
%
%
\newcount\subheadno%
\global\subheadno=0%
\def\subheadinfo{\the\subheadno}%
\def\nextsubheadno{\global\advance\subheadno by 1 \global\procno=0 \subheadinfo}%
\makecounter{subhead}%
%
%
\newcount\procno%
\global\procno=0%
\def\procinfo{\subheadinfo.\the\procno}%
\def\nextprocno{\global\advance\procno by 1 \procinfo}%
\makecounter{proc}%
%
%
\newcount\figno%
\global\figno=0%
\def\figinfo{\subheadinfo.\the\figno}%
\def\nextfigno{\global\advance\figno by 1 \figinfo}%
\makecounter{fig}%
%
%
\newcount\eqnno%
\global\eqnno=0%
\def\eqninfo{\text{{\rm (\the\eqnno)}}}%
\def\nexteqnno[#1]{\global\advance\eqnno by 1 \eqninfo\hbox{\eqnlabel{#1}}}%
\makecounter{eqn}%
%
%
%
%
%
\def\gobbleeight#1#2#3#4#5#6#7#8{}%
\newcount\citationno%
\global\citationno=0%
\def\citationinfo{\the\citationno}%
\makecounter{citation}%
\newwrite\biblio%
\def\newref#1#2{%
\def\temptext{#2}%
\edef\bibliotextoutput{\expandafter\gobbleeight\meaning\temptext}%
\global\advance\citationno by 1\citationlabel{#1}%
\ifmakebiblio%
    \edef\fileoutput{\write\biblio{\noindent\hbox to 0pt{\hss$[\the\citationno]$}\hskip 0.2em\bibliotextoutput\medskip}}%
    \fileoutput%
\fi}%
\def\cite#1{%
$[\citationref{#1}]$%
\ifmakebiblio%
    \edef\fileoutput{\write\biblio{#1}}%
    \fileoutput%
\fi%
}%
%
%
%
%
\let\mypar=\par%
\edef\Pagetitle={Blank}\headline={\hfil\Pagetitle\hfil}%
\edef\Pagefooter={Blank}\footline={\hfil\Pagefooter\hfil}%
%
%
\newcount\showpagenumflag%
\global\showpagenumflag=0 %
\def\nextoddpage%
{\newpage\ifodd\pageno%
\else\global\showpagenumflag=0 %
\null\vfil\eject%
\global\showpagenumflag=1 %
\fi}%
%
%
\font\headfont=cmb12%
\def\newhead#1[#2]%
{\ifhmode\mypar\fi%
\ifnum\headno=0 \else\goodbreak\bigskip\fi%
{\headfont\noindent\nextheadno\ - #1.}\headlabel{#2}%
\nobreak\medskip}%
%
%
\def\newsubhead#1[#2]%
{\ifhmode\mypar\fi%
\ifnum\subheadno=0 \else\goodbreak\medskip\fi%
{\bf\noindent\nextsubheadno\ - #1.\ }\subheadlabel{#2}}%
%
%
\newif\ifinproclaim%
\global\inproclaimfalse%
\def\proclaim#1{%
\goodbreak\medskip
\bgroup\inproclaimtrue%
\noindent{\bf #1}%
\nobreak\medskip\sl}%
\def\noskipproclaim#1{%
\goodbreak\medskip%
\bgroup\inproclaimtrue%
\noindent{\bf #1}\nobreak\sl}%
\def\endproclaim{\mypar\egroup\nobreak\medskip\ignorespaces}%
%
%
%
\newcount\xpos\newcount\ypos
\def\makelabelgrid{%
\xpos=-5 \ypos=-5 %
\loop\ifnum\xpos<6 %
{\loop\ifnum\ypos<6 %
\def\labeltext{x}%
\ifnum\xpos=0\def\labeltext{+}\fi%
\ifnum\ypos=0\def\labeltext{+}\fi%
\placelabel[\xpos][\ypos]{\labeltext}%
\advance\ypos by 1 %
\repeat}%
\advance\xpos by 1 %
\repeat}%
\def\placelabel[#1][#2]#3{{%
\setbox10=\hbox{\raise #2cm \hbox{\hskip #1cm #3}}%
\ht10=0pt \dp10=0pt \wd10=0pt \box10}}%
\def\placefigure#1#2#3#4{%
\medskip%
\midinsert%
\epsfxsize=3.5cm%
\vbox{\line{\hfil#2\epsfbox{#3}#1\hfil}%
\vskip 0.3cm%
\line{\hfil\vbox{\hsize=10cm \noindent{\sl Figure \nextfigno\ - #4}\par}\hfil}}%
\medskip%
\endinsert}%
%
%
\def\myitem#1{\noindent\hbox to .5cm{\hfill#1\hss}}%
%
%
%
%
%
%
%
%
%
\font\sansseriften=cmss10%
\font\sansserifseven=cmss7%
\font\sansseriffive=cmss5%
\newfam\sansseriffam%
\textfont\sansseriffam=\sansseriften%
\scriptfont\sansseriffam=\sansserifseven%
\scriptscriptfont\sansseriffam=\sansseriffive%
\def\mathsf{\fam\sansseriffam}%
%
%
%
\font\boldten=cmb10%
\font\boldseven=cmb7%
\font\boldfive=cmb5%
\newfam\mathboldfam%
\textfont\mathboldfam=\boldten%
\scriptfont\mathboldfam=\boldseven%
\scriptscriptfont\mathboldfam=\boldfive%
\def\mathbf{\fam\mathboldfam}%
%
%
%
\font\mycmmiten=cmmi10%
\font\mycmmiseven=cmmi7%
\font\mycmmifive=cmmi5%
\newfam\mycmmifam%
\textfont\mycmmifam=\mycmmiten%
\scriptfont\mycmmifam=\mycmmiseven%
\scriptscriptfont\mycmmifam=\mycmmifive%
\def\hexa#1{\ifcase #1 0\or 1\or 2\or 3\or 4\or 5\or 6\or 7\or 8\or 9\or A\or B\or C\or D\or E\or F\fi}%
\mathchardef\mathi="7\hexa\mycmmifam7B%
\mathchardef\mathj="7\hexa\mycmmifam7C%
%
%
\font\mymsbmten=msbm10 at 8pt%
\font\mymsbmseven=msbm7 at 5.6pt
\font\mymsbmfive=msbm5 at 4pt%
\newfam\mymsbmfam%
\textfont\mymsbmfam=\mymsbmten%
\scriptfont\mymsbmfam=\mymsbmseven%
\scriptscriptfont\mymsbmfam=\mymsbmfive%
\mathchardef\mybeth="7\hexa\mymsbmfam69%
\mathchardef\mygimmel="7\hexa\mymsbmfam6A%
\mathchardef\mydaleth="7\hexa\mymsbmfam6B%
%
%
%
%
\def\proof{{\noindent\bf Proof:\ }}%
\def\qed{~$\square$}%
\def\makeop#1{\global\expandafter\def\csname op#1\endcsname{{\text{#1}}}}%
\def\makeopsmall#1{\global\expandafter\def\csname op#1\endcsname{{\text{\lowercase{#1}}}}}%
%
%
%
%
%
%
\makeop{Ext}%
\makeop{Int}%
\makeop{Dist}%
\makeop{Diam}%
\makeop{Length}%
%
%
%
%
%
\def\mlim{\mathop{{\text{Lim}}}}%
%
%
%
%
%
%
\makeop{Dim}%
\makeop{Ker}%
\makeop{Coker}%
\makeop{Tr}%
\makeop{Adj}%
\makeop{Det}%
\makeop{End}%
\makeop{Lin}%
\makeop{Symm}%
\makeop{Mult}%
%
%
\makeop{dx}%
\makeop{dy}%
\makeop{dz}%
\makeop{dt}%
\makeop{dVol}%
\makeop{dArea}%
\makeop{Supp}%
\makeop{Hess}%
\makeop{Lip}%
%
%
\makeop{Re}%
\makeop{Im}%
\makeop{Arg}%
\makeop{Log}%
\makeop{Exp}%
%
%
\makeopsmall{Cos}%
\makeopsmall{Sin}%
\makeopsmall{Tan}%
\makeopsmall{Sec}%
\makeopsmall{Cosec}%
\makeopsmall{Cot}%
\makeopsmall{ArcCos}%
\makeopsmall{ArcSin}%
\makeopsmall{ArcTan}%
\makeopsmall{ArcSec}%
\makeopsmall{ArcCosec}%
\makeopsmall{ArcCot}%
%
%
\makeopsmall{Cosh}%
\makeopsmall{Sinh}%
\makeopsmall{Tanh}%
\makeopsmall{ArcCosh}%
\makeopsmall{ArcSinh}%
\makeopsmall{ArcTanh}%
%
%
\makeop{Vol}%
\makeop{Area}%
\makeop{Riem}%
\makeop{Ric}%
\makeop{Scal}%
\makeop{Euc}%
\makeop{Imm}%
\makeop{Emb}%
%
%
\makeop{Id}%
\makeop{Ad}%
\makeop{O}%
\makeop{SO}%
\makeop{SL}%
\makeop{GL}%
\makeop{Conf}%
\makeop{Homeo}%
\makeop{Diff}%
\makeop{Isom}%
%
%
\makeop{Ind}%
\makeop{Sig}%
\makeop{Spec}%
%
%
\makeop{Conv}%
\makeop{Max}%
\makeop{Min}%
\makeop{Mod}%
\makeop{Deg}%
\makeop{loc}%
%
%
%
%
%
%
%
%
%
%
%
%
%

%% file: references.tex
\global\def\_@citation@Brendle{1}
\global\def\_@citation@ChenFraserPang{2}
\global\def\_@citation@CodaNeves{3}
\global\def\_@citation@Devyver{4}
\global\def\_@citation@FraserSchoenI{5}
\global\def\_@citation@FraserSchoenII{6}
\global\def\_@citation@MaximoNunesSmith{7}
\global\def\_@citation@Hung{8}
\global\def\_@citation@Urbano{9}
\global\def\_@subhead@MorseIndexOfTheCriticalCatenoid{1}
\global\def\_@eqn@RobinBoundaryCondition{\relax \unhbox \voidb@x \hbox {{\relax \tenrm (1)}}}
\global\def\_@proc@MorseIndex{1.1}
\global\def\_@subhead@ProofOfTheoremX{2}
\global\def\_@eqn@EigenfunctionEquation{\relax \unhbox \voidb@x \hbox {{\relax \tenrm (2)}}}
\global\def\_@eqn@OneDimRobin{\relax \unhbox \voidb@x \hbox {{\relax \tenrm (3)}}}
\global\def\_@eqn@EEInFourierMode{\relax \unhbox \voidb@x \hbox {{\relax \tenrm (4)}}}
\global\def\_@eqn@RCInFourierMode{\relax \unhbox \voidb@x \hbox {{\relax \tenrm (5)}}}
\global\def\_@proc@UpperLimits{2.1}
\global\def\_@eqn@RicattiEquation{\relax \unhbox \voidb@x \hbox {{\relax \tenrm (6)}}}
\global\def\_@fig@FigureA{2.1}
\global\def\_@subhead@Bibliography{3}